\documentclass{elsart}

\usepackage{natbib}

\usepackage{amssymb}

\newtheorem{lemma}[thm]{Lemma}
\newtheorem{corol}[thm]{Corollary}
\newtheorem{propos}[thm]{Proposition}

\def\reff#1{(\ref{#1})}
\def\one{{\mathbf 1}}
\def\emph#1{{\it#1}}
\def\p{q}

\def\wp{\widetilde p}

\def\wL{\widetilde L}
\def\wH{\widetilde H}
\def\wmu{\widetilde \mu}

\def\wp{\widetilde p}
\def\wP{\widetilde P}

\def\E{{\mathbb E}}
\def\P{{\mathbb P}}
\def\R{{\mathbb R}}
\def\Z{{\mathbb Z}}

\def\L{\Lambda}

\def\eps{\varepsilon}

\def\1{{1\kern-.25em\hbox{\rm I}}}
\def\eu{{1\kern-.25em\hbox{\sm I}}}

\def\R{{\mathbb R}}  
\def\P{{\mathbb P}}  
\def\Z{{\mathbb Z}}  
\def\E{{\mathbb E}}  

\let\cal=\mathcal

\def\EE{{\cal E}}
\def\FF{{\cal F}}

\def\II{{\cal I}}
\def\HH{{\cal H}}

\def\TT{{\cal T}}

\def\UU{{\cal U}}

\def\XX{{\cal X}}

\def\part#1{
  \shipout\vbox{
    \vbox to 0.35\vsize{}
    \hugetitle #1
    }
  \advancepageno
  \blankpage
}

\def\text#1{\quad{\hbox{#1}}\quad}

\def\proof{{\vskip 2mm \noindent\bf Proof }}

\def\sqr{\vcenter{
         \hrule height.1mm
         \hbox{\vrule width.1mm height2.2mm\kern2.18mm\vrule width.1mm}
         \hrule height.1mm}}                  

\def\thanks{\noindent{\bf Acknowledgements: }}

\def\proof{{\noindent\bf Proof. }}

\begin{document}

\begin{frontmatter}

\title{Harness processes and harmonic crystals}

\author{Pablo A. Ferrari},
\author{Beat M. Niederhauser}
\address{Universidade de S\~{a}o Paulo}

\begin{abstract}
In the Hammersley harness processes the $\R$-valued height at each
site $i\in\Z^d$ is updated at rate 1 to an average of the neighboring
heights plus a centered random variable (the noise). We construct the
process ``a la Harris" simultaneously for all times and boxes
contained in $\Z^d$.  With this representation we compute covariances
and show $L^2$ and almost sure time and space convergence of the
process. In particular, the process started from the flat
configuration and viewed from the height at the origin converges to an
invariant measure. In dimension three and higher, the process itself
converges to an invariant measure in $L^2$ at speed $t^{1-d/2}$ (this
extends the convergence established by Hsiao). When the noise is
Gaussian the limiting measures are Gaussian fields (harmonic crystals)
and are also reversible for the process.  \bigskip
\end{abstract}

\end{frontmatter}

\section{Introduction}

\paragraph*{\bf The harness process} 
The {\sl harness process} is a continuous-time version of the serial harness
introduced by Hammersley \cite{H}. Let $P=(p(i,j), i,j\in\Z^d)$ be a translation
invariant finite-range stochastic matrix (that is, $p(i,j)\ge 0$,
$\sum_jp(i,j)=1$ for all $i$, $p(i,i+j)= 0$ if $|j| > v$ for some $v$ and
$p(i,j) = p(0,j-i)$ for all $i,j$). Let the \emph{noise} $G(dx)$ be a centered
distribution with variance $1$. The state-space is $\XX=\R^{\Z^d}$. We consider
a family of processes in subsets $\Lambda\subset\Z^d$ with boundary conditions
$\gamma\in\XX$. For configurations $\eta\in\XX$ and bounded cylinder functions
$f:\XX\to\R$ define the generator
\begin{equation}
  \label{p100}
  L^{\Lambda,\gamma} f(\eta) = \sum_{i\in\Lambda} \int G(d\eps)
  [f(P_i(\eta_\Lambda\gamma_{\Lambda^c})+\sigma\eps e_i) - f(\eta)] 
\end{equation}
where the standard deviation of the noise $\sigma>0$ is a parameter, $e_i(j) =
\1\{i=j\}$, $P_i\eta$ is the configuration
\begin{equation}
  \label{p2}
  (P_i\eta)(i)  = \sum_{j\in\Z^d} p(i,j) \eta(j)\,;\quad (P_i\eta)(j) =
  \eta(j) \hbox{ for } j\neq i\,;\qquad
\end{equation}
and the juxtaposition $\eta_\Lambda \gamma_{\Lambda^c}\in\XX$ is defined by
\begin{equation}
  \label{p70}
  (\eta_A\gamma_{A^c})(i) \,=\,\left\{\begin{array}{ll}
\eta(i), & \hbox{if\ }i\in A, \cr
       \gamma(i), & \hbox{if\ }i\in A^c. \cr
\end{array}\right.
\end{equation}
In other words, at all times the sites outside $\Lambda$ have fixed
configuration $\gamma$ and those inside are updated at rate 1 with a
$P$-weighted mean of the neighbors plus an independent centered random variable.
When the boundary configuration $\gamma$ is the \emph{flat configuration}
$\gamma(i)\equiv 0$ we write $L^\Lambda$.

Basis \cite{Ba1} \cite{Ba2} proves that there exist a Markov processes
$(\eta_t)$ in $\R^{\Z^d}$ with generators $L^{\Lambda,\gamma}$, that
is, processes satisfying
\begin{equation}
  \label{p4}
  \lim_{h\to 0} \frac{1}{h}\E [f(\eta_{t+h}) -
  f(\eta_t)\,|\,\FF_t]  
  = L^{\Lambda,\gamma}f(\eta_t)
\end{equation}
for bounded cylinder functions $f$, where $\FF_t$ is the
$\sigma$-algebra generated by $\{\eta_s,\,s\le t\}$. His proof works
in a more general context of metric spaces. The existence is
immediate if $\Lambda$ is finite but for infinite $\Lambda$ it is
necessary to impose the boundary conditions $\gamma$ not to grow too
fast (see \reff{o1} later). 
Hsiao \cite{Hs1} \cite{Hs2} shows existence of invariant
measures in dimensions $d\ge 3$ and gives conditions for the
convergence of the process to the invariant measures.  The
discrete-time version is called ``serial-harness'' by Hammersley and
its tail behavior has been studied by Toom~\cite{T}.

\paragraph*{\bf The Gaussian Gibbs fields}
For each finite $\Lambda\subset \Z^d$ let $H^\Lambda:\XX\to \R$ be the Hamiltonian
\begin{equation}
  \label{I.14.1}
 H^\Lambda(\eta) =  \frac{\beta}{2} \sum_{i\in\Lambda}\sum_{j\in\Z^d} 
p(i,j)(\eta(i)- \eta(j))^2 
\end{equation}

For finite  $\Lambda\subset\Z^d$ and $\gamma\in\XX$ define the
measure $\mu^{\Lambda,\gamma}$ on $\R^\Lambda$ by
\begin{equation}
  \label{p71}
  \mu^{\Lambda,\gamma} (f) = \frac{1}{Z^{\Lambda,\gamma}}
  \int_{\R^\Lambda} f(\eta)e^{- H^\Lambda(\eta_\Lambda\gamma_{\Lambda^c})} 
  \prod_{i\in \Lambda}d\eta(i)
\end{equation}
where $d\eta(i)$ is the Lebesgue measure in the $i$th coordinate of
$\Lambda$. The elements of the family
\begin{equation}
  \label{p72}
  \{\mu^{\Lambda,\gamma}\,:\, \Lambda\subset\Z^d \hbox{ finite},\; \gamma\in\XX\}
\end{equation}
are called \emph{local specifications}. When $\gamma$ is the flat
configuration we write $\mu^\Lambda$.  One of the main problems in
Statistical Mechanics is to find a measure on $\XX$ whose conditional
probabilities are given by the specifications \reff{p71} (DLR
equations, see the book of Georgii \cite{G} or the monograph of
Bovier \cite{B}; for the Gaussian fields this has been solved by
Spitzer \cite{Sp2} and Dobrushin \cite{D}).  More precisely, we say
that a measure $\mu$ is a Gibbs measure with specifications $
\mu_{\Lambda,\gamma}$ if for all finite $\Lambda$ and continuous
$f:\R^\Lambda\to\R$, the conditional probabilities exist $\mu$ almost
surely and satisfy
\begin{equation}
  \label{p73}
  \mu(\cdot|\FF_\Lambda^c)(\gamma_{\Lambda^c}) = \mu^{\Lambda,\gamma}
  \qquad\mu\hbox{ a.s.}
\end{equation}
where $\FF_\Lambda^c$ is the $\sigma$-algebra generated by $\gamma_{\Lambda^c}$.

\paragraph*{\bf Harnesses}
The motivation of Hammersley \cite{H} was the construction of probability
measures $\mu$ on $\R^{\Z^d}$ with the property 
\begin{equation}
  \label{h1}
  \mu(\eta(x)\,|\,\eta(y),\,y\neq x) = \sum_y p(x,y) \eta(y) 
\end{equation}
that is, the expected value under $\mu$ of the height at $x$ conditioned on
the heights at the other sites is a convex combination (taken with the matrix
$p$) of the heights at the other sites. Measures $\mu$ satisfying \reff{h1}
are called \emph{harnesses}. Williams \cite{W} constructs Gaussian measures
that are harnesses when $p$ is a nearest neighbor symmetric random walk in
$\Z^d$. Kingman \cite{K} proposes the construction of harnesses in $L^1$.
The Gaussian Gibbs fields satisfying \reff{p73} are harnesses.

\paragraph*{\bf Results}
The point of this paper is a simultaneous construction (coupling) of versions of
the processes $(\eta^{\Lambda,\gamma}_t)$ and configurations
$\eta^{\Lambda,\gamma}$ with law $\mu^{\Lambda,\gamma}$ for all $\Lambda$ and
$\gamma$, in the same probability space. Then we show $L_2$ and almost sure time
and space convergence. This is based on a Harris graphical construction of the
harness process on a probability space generated by a family of one-dimensional
marked stationary Poisson processes indexed by $\Z^d$. Epochs of the Poisson
process correspond to updating times of the Harness process; the marks are
independent and identically distributed random variables with distribution $G$.
This construction allows to represent the process starting at time $s$ with the
flat configuration as
\begin{equation}
  \label{p20i}
  \eta_{[s,t]}(i) \,:=\, \sum_{j\in\Lambda}\;
  \sum_{n:T_n(j)\in [s,t]} \eps_n(j)\,b_n(i,j)
\end{equation}
for $t\ge s$. Here $\eps_n(j)$ is the noise associated to $T_n(j)$, the $n$th
Poisson epoch of site $j$ and $b_n(i,j)$ is the probability that \emph{given the
Poisson epochs}, a random walk starting at time $t$ at site $i$ jumping at the
Poisson epochs backwards in time is at site $j$ at time $T_n(j)$. The jumps of
the walk have law $p$. Since $b_n(i,j)$ are a function of the Poisson epochs,
$\eta_{[s,t]}(i)$ is a function of the Poisson epochs in the interval $[s,t]$
and the noises associated to them.  This representation is the continuous
analogous of equation (8.2) in \cite{H}. It is reminiscent of what is called
\emph{duality} in interacting particle systems and goes in parallel with the
backwards representation of the \emph{random average process} in \cite{FF}.

We show that for each fixed $t$ the process $(\eta_{[t-s,t]}(i)\,,\,s\ge 0)$ is
a martingale with uniformly bounded second moments in $d\ge 3$ and hence for
each fixed $t$ it converges almost surely to a limit denoted $\eta_t(i)$. We
also show that the rate of $L^2$ convergence is bounded by a constant times
$s^{1-d/2}$, improving the weakly convergence established by Hsiao \cite{Hs1}
\cite{Hs2}. The limiting process $(\eta_t,\,t\in\R)$ is a stationary harness
process.  In $d\le 2$ we study the process pinned at zero in the origin (for
which the site at the origin is not updated and remains zero) and the process as
seen from the height at the origin. We prove similar results in those cases. To
our knowledge these results are new in $d=1,2$. The graphical construction and
the martingale property are shown in Section \ref{s2}.

The process can be defined in subsets of $\Lambda\subset\Z^d$ by
assuming that the heights outside $\Lambda$ are fixed. Using the
superlabel $\Lambda$ for the process restricted to $\Lambda$ with the
heights outside $\Lambda$ equal to zero we get a family of stationary
processes $((\eta^\Lambda_t,\,t\in\R),\,\Lambda\subset\Z^d)$. We show
that under suitable conditions, for each $t$, the one-time marginal
family $(\eta^\Lambda_t,\,\Lambda\subset\R^d)$ converges coordinatewise
in $L^2$ to an infinite volume configuration $\eta^{\Z^d}_t$ as
$\Lambda\nearrow\Z^d$.  The time and space convergence results are
proven in Theorem \ref{x1} in Section~\ref{s4}.

The one and two-point correlations are computed in Section \ref{s3} using the
following random walk representation of the second moments of the differences:
For $i\in\Z^d$,
\begin{equation}
  \label{F.19x}
  \E\,(\eta_{[s,t]}(i) - \eta_{[s,t]}(0)) ^2\, = \,
  2\,\int_0^{t-s} (\P(D^0_u= 0) - \P(D^{i}_u= 0))\,du
\end{equation}
where $D^i_u$ is the position at time $u$ of a symmetric random walk starting at
$i$ at time $0$. The transition probabilities of this walk are homogeneous but
at the origin; they are given in \reff{F.17}. This walk also appears in Hsiao
\cite{Hs1} to compute the correlations of the stationary law of $\eta_t$.

The law of $\eta^\Lambda_t$ is the unique invariant measure for the
harness process when $\Lambda$ is finite; recall that the boundary
conditions we are taking ``pin'' the process to the external
configuration. This is proven in Theorem \ref{x1} using the
representation \reff{p20i}. In the infinite case there are infinitely
many invariant measures. In particular, if $h$ is a harmonic function
for $p$, in $d\ge 3$ the law of $\eta^{\Z^d}_t+h$ is invariant for the
harness process.  We conjecture that in $d\ge 3$ the law of
$\eta^{\Z^d}_t$ is the unique ergodic invariant measure with mean
zero. Hsiao \cite{Hs1} proved that this is the only ergodic invariant
measure with mean zero and finite variance.  To eliminate the
restriction of finite variance it would be sufficient to show the
following random version of the ergodic theorem: Let $\eta$ be a
configuration chosen from an ergodic measure $\mu$ with mean zero and
$\P$ the probability induced by the Poisson processes, then
\begin{equation}
  \label{z14}
   \lim_{s\to\infty} \sum_j b_{[0,-s]}(i,j) \eta(j) = 0 \qquad{\rm
   \P\hbox{-}a.s.\; 
  \mu\hbox{-}a.s. }
\end{equation}
where $b_{[0,-s]}(i,j)$ is the probability conditioned on the Poisson epochs
that the backwards walk starting at $i$ at time 0 is at $j$ at time $-s$. The
ergodicity of $\mu$ implies that \reff{z14} holds $\mu$-a.s.\/ if we replace
$b_{[0,-s]}(i,j)$ by its averages. The limit \reff{z14} is related to the
asymptotic behavior of the \emph{no-noise} harness process $\underline\eta_t$
defined in \reff{x35}, a harness process with zero noise (that is,
$G(dx)=\delta_0(x)$). In this process the heights are updated at the Poisson
times to the $p$-average of the other heights. The problem is to characterize
the set of initial configurations for which this process converges to the
``all-zero'' configuration.

Under the assumptions that the noise $G$ is Gaussian (that is $G(dx) =
(2\pi)^{-1/2}$ $e^{-x^2/2} dx$) and that $p(0,0)=0$, Hsiao \cite{Hs1} proved
that the Gaussian Gibbs field $\mu^\Lambda$ is reversible for the harness
process in any $\Lambda$. Indeed, since the conditional distribution under
$\mu^\Lambda$ of $\eta(i)$ given $(\eta(j),\,j\neq i)$ has Gaussian law
centered at $\sum_j p(i,j)\eta(j)$, the harness process is just the so called
\emph{heat bath} dynamics at continuous time. The weak convergence of
$\mu^\Lambda$ to $\mu^{\Z^d}$ for $d\ge 3$ has been proven by Spitzer
\cite{Sp2}; we provide here convergence in $L^2$ and a simultaneous
construction of $(\xi^\Lambda)_\Lambda$ for an increasing sequence of finite
sets $\Lambda\nearrow\Z^d$ satisfying that $\xi^\Lambda$ has law $\mu^\Lambda$
and converges almost surely to a configuration $\xi^{\Z^d}$ with law
$\mu^{\Z^d}$, the infinite volume Gibbs measure with specifications
\reff{p71}. This is done in Proposition \ref{xx4}. The almost sure convergence of
$\eta^\Lambda_t$ as $\Lambda\nearrow\Z^d$ in $d\ge 3$ remains open. We prove
similar results for the process pinned at the origin and the process as seen
from the height at the origin.

Compared with the work of Hsiao who considered $d\ge 3$, our constructive
approach permits (a) to treat (bounded or unbounded) regions $\Lambda$ contained
in $\Z^d$ and the difference process in dimensions $d=1$ and $2$ and (b) compute
non-equilibrium correlation functions.  Hsiao also considered the case when $p$
is sub-stochastic; we discuss this with Pechersky \cite{FNP}.

\section{Harris graphical construction}\label{s2}
Let $(\TT,\EE,\UU)$ be a collection of independent marked rate-1 Poisson
processes on $\R$:
\begin{eqnarray}
  \label{o4}
   (\TT,\EE, \UU) :=((T_n(i),\eps_n(i),U_n(i))\,;\; i\in \Z^d,\,n\in \Z)
\end{eqnarray}
where $T_n(i)$ is the $n$th epoch of a stationary Poisson process of rate $1$
(that is, $T_0(i)<0\le T_1(i)$, $T_1(i)$, $-T_0(i)$ and $T_n(i)-T_{n-1}(i)$ for
$n\neq 1$ are i.i.d.\/ exponential with mean 1); $\eps_n(i)$ are i.i.d.\/
centered random variables with variance $1$ and $U_n(i)$ are i.i.d.\/ in $\Z^d$
with law $p(i,\cdot)$. Furthermore $T_n(i)-T_{n-1}(i)$, $\eps_{n'}(i')$,
$U_{n''}(i'')$, $n,n',n'',i,i',i''\in \Z$ are mutually independent random
variables. Let $\P$ and $\E$ denote the probability and expectation induced by
these processes.

Fix $t\in\R$ and let $(B^{i,\Lambda}_{[t,u]}\,,\,u\le t)$ be a \emph{backward}
random walk starting at site $i$ at time $t$ and jumping at the Poisson epochs
backwards in time according to the $U_n(j)$ variables and \emph{absorbed} at
$\Lambda^c$. That is, $B^{i,\Lambda}_{[t,t]}=i$ and if at time $u+$ the walk
is at $j\in\Lambda$, $T_n(j)=u$ and $U_n(j)=j'$, then at this time the walk
jumps to $j'$. If $j\notin\Lambda$ then it stays at $j$ for ever.

For $s\le t$ define $\eta^\Lambda_{[s,t]}(i)$ as the expectation of the sum of
the noise variables $\varepsilon_n(i)$ encountered by
$B^{i,\Lambda}_{[t,\cdot]}$ in the (backwards) interval $[t,s]$
\emph{conditioned on the jump times}.  More precisely, define
$\eta^\Lambda_{[s,s]}(i)\equiv 0$ and for $t\ge s$,
\begin{equation}
  \label{p20}
  \eta^\Lambda_{[s,t]}(i) \,:=\, \sum_{j\in\Lambda}\;
  \sum_{n:T_n(j)\in [s,t]} \eps_n(j)\,b^\Lambda_{[t,T_n(j)]}(i,j)
\end{equation}
where, abusing notation by calling $\TT$ the $\sigma$-algebra generated by
$\TT$, 
\begin{equation}
  \label{p21}
  b^{\Lambda}_{[t,u]}(i,j)\,=\,b^{\Lambda}_{[t,u]}(i,j\vert\TT)\,:=\,
  \P(B^{i,\Lambda}_{[t,u]} = j\,|\, \TT) 
\end{equation}
for $u\le t$; that is, $b^{\Lambda}_{[t,u]}(i,j)$ is a function of the Poisson
epochs in the interval $[u,t]$ and it is independent of $\EE$ and $\UU$. 

For each $s\in\R$, expressions \reff{p20} and \reff{p21} define a random process
$(\eta^\Lambda_{[s,t]},\, t\ge s)$ as a (deterministic) function of
$((T_n(j),\eps_n(j))\,:\, T_n(j)\in [s,\infty))$. The sums \reff{p20} are almost
surely finite as a consequence of the finite range of $p$ and the fact that
there are only a finite number of Poisson epochs in bounded time intervals.

We also define the process starting with a configuration $\zeta\in\XX$ at time
$s$ by $\eta^{\Lambda,\zeta}_{[s,s]}(i)\equiv \zeta$ and for $t\ge s$,
\begin{equation}
  \label{p20p}
  \eta^{\Lambda,\zeta}_{[s,t]}(i) \,:=\, \sum_{j\in\Lambda}\;
  \sum_{n:T_n(j)\in [s,t]} \eps_n(j)\,b^\Lambda_{[t,T_n(j)]}(i,j)\, +
  \,\sum_{j\in\Lambda} b^\Lambda_{[t,s]}(i,j)\,\zeta(j)
\end{equation}
This is defined for configurations $\zeta$ that do not increase too fast
to guarantee that the sum in \reff{p20p} is almost sure finite. A
sufficient condition is that $\zeta$ belongs to $\Xi_\Lambda$, where 
\begin{equation}
  \label{o3}
  \Xi_\Lambda :=\{\zeta\,:\,\sum_{j\in\Lambda}
  p^\Lambda_t(i,j)\,\zeta(j)<\infty\hbox{ for all }i\in\Lambda,\,t>0\}.
\end{equation}
where $p^\Lambda_t$ is the probability that a continuous random walk
with rates $p$, absorbed at sites in $\Lambda^c$ starting at $i$ at
time zero is at $j$ at time $t$. Notice that $p^\Lambda_{t-s} =
\E(b^\Lambda_{[t,s]}(i,j))$.

\begin{propos}
  \label{xx2}
For any $d\ge 1$, $\Lambda\subset\Z^d$ and $s\in\R$, the process
$(\eta^{\Lambda,\zeta}_{[s,t]},\, t\ge s)$ defined in \reff{p20} has generator
$L^{\Lambda}$ (in the sense of \reff{p4}) and initial condition $\zeta$ at time
$s$.
\end{propos}

\proof For any $s\in\R$ the process $(\eta^{\Lambda,\zeta}_{[s,t]},\, t\ge s)$
as defined by \reff{p20p} satisfies the following infinitesimal evolution:
\begin{equation}
  \label{I.13}
  \eta^{\Lambda,\zeta}_{[s,t]}(i)= \left\{\begin{array}{ll}  \eta^{\Lambda,\zeta}_{[s,t-]}(i), 
    &\hbox{ if $t$ is not 
          a epoch of $\TT(i)$}\\
      \sum_{j \in \Z^d} p(i,j) \eta^{\Lambda,\zeta}_{[s,t-]}(j) + \varepsilon_n(i),
                                & \hbox{ if $t = T_n(i)$}
         \end{array}\right.
\end{equation}
from where it follows that $\eta^{\Lambda,\zeta}_{[s,t]}$ has generator $L$.
$\sqr$

In the sequel we use the notation:
\begin{equation}
  \label{x31}
  b^{\Lambda}_{n}(i,j)\;:=\;b^{\Lambda}_{[t,T_n(j)]}(i,j)
\end{equation}

\begin{propos}\label{p777}
  For each $i\in\Z^d$ and $t\in\R$ the process
  $(\eta^\Lambda_{[t-s,t]}(i)\,,\,s\ge 0)$ is a martingale with respect to the
  filtration $(\FF_s)_{s \geq 0}$, where $\FF_s$ is the sigma algebra
  generated by $(\eta^\Lambda_{[t-u,t]})_{u \leq s}$.
\end{propos}

\proof For $r>s$ the expectation of $\eta^\Lambda_{[t-r,t]} -
\eta^\Lambda_{[t-s,t]} $ given $\FF_s$ vanishes because it is the mean of a
(random) finite sum of randomly weighted centered variables $\eps_n(j)$
independent of the weights and of the past.  Indeed, for $0 \le s\le r$,
\begin{eqnarray}
  \label{p28}
\lefteqn{  \E(\eta^\Lambda_{[t-r,t]} -
\eta^\Lambda_{[t-s,t]}\,|\,\FF_s)}\nonumber\\ 
  &=&  \E\biggl[\sum_{j\in\Z^d}\;\sum_{n:T_n(j)\in [t-r,t-s]}
  \eps_n(j)\,b^\Lambda_n(i,j)\bigg| \FF_s\biggr]\nonumber \\
  &=&  \E\biggl[\E\biggl[\sum_{j\in\Z^d}\;\sum_{n:T_n(j)\in [t-r,t-s]}
  \eps_n(j)\,b^\Lambda_n(i,j)\bigg|\TT,\FF_s\biggr]\bigg| \FF_s\biggr]\nonumber \\
  &=&  \E\biggl[\sum_{j\in\Z^d}\;\sum_n 
  \E[\eps_n(j)\,|\, \TT,\FF_s]\,b^\Lambda_n(i,j)\,\1\{n:T_n(j)\in [t-r,t-s]\}
  \bigg| \FF_s\biggr]\;=\;0 \nonumber
\end{eqnarray}
where the third identity follows from Fubini and the fact that both
$\p_n(i,j)$ and $T_n(j)$ are $\TT$-measurable; the fourth identity follows
because (a) for $T_n(j)\in[t-u,t-s]$, $\eps_n(j)$ is independent of $\FF_s$,
(b) $\eps_n(j)$ is independent of $\TT$ for all $n$ and $j$ and (c)
$\eps_n(j)$ are centered random variables.   $\sqr$

\section{Covariances}\label{s3}
This section collects bounds for the relevant covariances. The main tool is an
expression of the covariances of the process in $\Z^d$ as a function of the
potential kernel of a symmetric random walk. These covariances are bounds for
the covariances in the box $\Lambda\subset\Z^d$; this works in $d\ge 3$. As a
consequence the relevant variances are uniformly bounded in time and space. When
$p(0,0)=0$, we use results from the Gaussian case to bound the variances when
$d=1,2$ for the process ``pinned'' at the origin and for the process ``as seen
from the height at the origin''. The results are summarized in Corollary
\ref{xx11} later.

We start with an elementary computation.
\begin{lemma}
  \label{3}
  Let $\Lambda'\subset \Lambda\subset \Z^d$ and
  $s'\le s\le t$. For all $i\in\Lambda$, $j\in\Lambda'$ 
\begin{equation}
  \label{F.1900}
  \E\,[\eta^\Lambda_{[s,t]}(i) \eta^{\Lambda'}_{[s',t]}(j)] \;=\; 
        \E\Bigl[\sum_{k\in\Lambda}\sum_{n:T_n(k)\in [s,t]} 
  b^\Lambda_n(i,k)\,b^{\Lambda'}_n(j,k)\Bigr]\,.
\end{equation}
\end{lemma}

\proof Using the definition, conditioning on the Poisson marks and integrating
with respect to the disorder variables, the left hand side of \reff{F.1900}
equals
\begin{eqnarray}
  \label{F.20}
\lefteqn{\E\Bigl[\E\Bigl[\Bigl(\sum_{k\in\Z^d}\;
\sum_{n:T_n(k)\in [s,t]} \eps_n(k)\,b^{\Lambda}_n(i,k)\Bigr)\;
\Bigl(\sum_{k\in\Z^d}\sum_{n:T_n(k)\in [s',t]}
\eps_n(k)\,b^{\Lambda'}_n(j,k)\Bigr) 
\Big|\TT\Bigr]\Bigr]} \nonumber \\
&=& \E\Bigl[\sum_{k\in\Z^d}\sum_{k'\in\Z^d}
\sum_{n:T_n(k)\in [s,t]}
\sum_{n':T_{n'}(k')\in [s',t]} 
  \E(\eps_n(k)\eps_{n'}(k')\vert \TT)
\,b^{\Lambda}_n(i,k)
\,b^{\Lambda'}_{n'}(j,k')\Bigr]\nonumber \\ 
\end{eqnarray}
where we can interchange sums and conditional expectations as the sums are
$\TT$ almost surely finite. From \reff{F.20} we get the right hand side of
\reff{F.1900} because $\eps_n(k)$ are i.i.d.\/ independent of $\TT$ with
variance 1. $\sqr$

\paragraph*{\bf Covariances in $\Z^d$}
Let $D^i_t$ be a continuous time random walk on $\Z^d$ starting at $i$ with
the following (symmetric) transition rates:
\begin{equation}
  \label{F.17}
  p_D(i,j) = \left\{\begin{array}{ll}  p(0, j - i)+p(0,i-j) 
    ,& \hbox{if $i \neq 0$;} \cr 
  \sum_{k \in \Z^d} p(0,k)p(0,k+j),
                  &\hbox{if $i = 0$.} \cr
                  \end{array}\right.
\end{equation}
\begin{lemma}
  \label{p77}
  Let
  $d\ge 1$ and $-\infty<s\le t$. For $i,j\in\Z^d$, $\E\eta^{\Z^d}_{[s,t]}(i)=0$
  and
\begin{eqnarray}
  &&  \E\,(\eta^{\Z^d}_{[s,t]}(i))^2 \;=\; 
        \E\Bigl[\sum_{k\in\Z^d}\sum_{n:T_n(k)\in [s,t]} 
  b^{\Z^d}_n(i,k)^2\Bigr] \;=\; 
         \int_0^{t-s} \,\P(D^0_u= 0)\,du\label{F.19a} \\
&&
  \E\,(\eta^{\Z^d}_{[s,t]}(j) - \eta^{\Z^d}_{[s,t]}(i)) ^2\, = \,
         2\,\int_0^{t-s} (\P(D^0_u= 0) - \P(D^{i-j}_u= 0))\,du
  \label{F.19}
\end{eqnarray}
\end{lemma}

\proof Taking $\Lambda = \Lambda' = \Z^d$, $s=s'$ in \reff{F.1900} gives the
first identity in \reff{F.19a}. The middle expression in \reff{F.19a} is the
average number of Poisson epochs used simultaneously by $B^{i,\Lambda}_{[t,s]}$
and $\bar B^{i,\Lambda}_{[t,s]}$, where $\bar B^{j,\Lambda}_{[t,s]}$ is a random
walk that uses the same Poisson epochs as $B^{i,\Lambda}_{[t,s]}$ but
independent jump variables $\bar U_n(\cdot)$.  Noting that $\bar
B^{j,\Lambda}_{[t,s]}-B^{i,\Lambda}_{[t,s]}$ has the same law as
$D^{i-j}_{t-s}$, the second identity in \reff{F.19a} follows. In this
computation the expected number of Poisson marks at the origin seen by
$D^{i-j}_u$, for $u\in[t,s]$ equals the right hand side of \reff{F.19a} because
the jump rate at the origin is 1. The same considerations show \reff{F.19}.
$\sqr$

We now get bounds for the time integrals.

\begin{lemma}
\label{4}
There exist constants $C$ and $C(i)$ such that for $s>1$,
\begin{eqnarray}
&&  \int_s^\infty \P(D^0_u= 0)\,du \;<\;C s^{1-d/2},\qquad \hbox{ for } d\ge 3\label{bb2}\\
&& \int_s^\infty (\P(D^0_u= 0) - \P(D^{i}_u= 0))\,du\;<\; C(i) s^{-d/2}\label{bb3}
\end{eqnarray}
\end{lemma}

\proof 
Since $D$ is a local perturbation of a symmetric finite range random walk, we
have $P(D^0_u= 0)<Cs^{-d/2}$, from where one gets \reff{bb2} with another constant. 
Differentiating with respect to $u$ and using Kolmogorov Backwards equation we
get 
\[
P(D^0_s=0)= \int_s^\infty \sum_i p_D(0,i) (\P(D^0_u= 0) - \P(D^{i}_u= 0)) \,du  
\]
Since the differences are positive, we get \reff{bb3} with $C(i)=C/p_D(0,i)$
when $p_D(0,i)>0$. An inductive step shows \reff{bb3} for all $i$.  $\sqr$

Next we show that if $p(0,0)=0$, the variances of the process pinned at zero are
uniformly bounded. The property holds for all centered noises of variance 1, but
the proof uses the fact that the Gibbs measure with specifications \reff{p71} is
reversible for the process with Gaussian noise. This is the case only when
$p(0,0)=0$. 
\begin{lemma}
  \label{x70}
  Assume $p(0,0)=0$. Then for all $d\ge 1$, $i\in\Lambda$,
  $\Lambda\subset\Z^d$ there exist constants
  $V^{\Lambda\setminus\{0\}}(i)<\infty$ such that
  \begin{equation}
    \label{x21x}
    \E\,[\eta^{\Lambda\setminus\{0\}}_{[s,t]}(i)]^2
\;\le\;V^{\Lambda\setminus\{0\}}(i)\;<\;\infty
  \end{equation}
\end{lemma}
\proof From \reff{F.1900} we see that the variances do not depend on the
particular distribution $G$ provided its variance is 1.  Hence we can assume
without loss of generality that the noise is Gaussian. Theorem \ref{xx3} later
  says that under $p(0,0)=0$ and Gaussian noise there exists a Gibbs measure
  $\mu^{\Lambda\setminus\{0\}}$ reversible (hence invariant) for the process.
  That is,
\begin{equation}
  \label{x73}
  \int\mu^{\Lambda\setminus\{0\}}(d\xi)
  \,\E f(\eta^{\Lambda\setminus\{0\},\xi}_{[s,t]})\; =
  \;\int\mu^{\Lambda\setminus\{0\}}(d\xi) f(\xi)
\end{equation}
for cylinder continuous $f:\XX\to\R$.  The variances
$V^{\Lambda\setminus\{0\}}(i)=:
\int\mu^{\Lambda\setminus\{0\}}(d\xi)\,\xi(i)^2$
are finite for all $\Lambda\subset\Z^d$ (see \reff{l12} later).
Then, using \reff{p20p} and the invariance property \reff{x73},
  \begin{eqnarray}
    \label{x71}
  V^{\Lambda\setminus\{0\}}(i)&=& \int\mu^{\Lambda\setminus\{0\}}(d\xi)
    \,\E\,[\eta^{\Lambda\setminus\{0\},\xi}_{[s,t]}(i)]^2\nonumber\\
 &=& 
    \E\,[\eta^{\Lambda\setminus\{0\}}_{[s,t]}(i)]^2 +
    \int\mu^{\Lambda\setminus\{0\}}(d\xi)\,\E\Bigl(\sum_{k\in\Lambda}
    b^{\Lambda\setminus\{0\}}_{[t,s]}(i,k)\,
     \xi(k)\Bigr)^2 
  \end{eqnarray}
  (The crossed terms cancel because $\eps_n(k)$ are centered and independent
  of $\xi$ and $b$.) This shows \reff{x21x}. $\sqr$

  Variances are monotone in time and $\Lambda$:

\begin{lemma}
  \label{x11}
For $i\in\Z^d$, $\Lambda\subset\bar\Lambda$ and $ t\ge s\ge \bar s$, 
\begin{equation}
  \label{F.19b}
  \E\,[\eta^{\Lambda}_{[s,t]}(i)]^2\, \le
  \,\E\,[\eta^{\bar\Lambda}_{[\bar s,t]}(i)]^2\, 
\end{equation}
  \begin{equation}
    \label{x20}
\E\,[\eta^{\Z^d}_{[s,t]}(i) -  \eta^{\Z^d}_{[s,t]}(0)]^2 
\;\le\;
\E\,[\eta^{\Z^d}_{[\bar s,t]}(i) -  \eta^{\Z^d}_{[\bar s,t]}(0)]^2 
  \end{equation}
\end{lemma}
\proof Using \reff{F.1900} with $\Lambda=\Lambda'$ and $s=s'$:
\begin{eqnarray}
  \label{p91}
  \E\,[\eta^{\Lambda}_{[s,t]}(i)]^2&=&
  \E\Bigl[\sum_{k\in\Lambda}\sum_{n:T_n(k)\in [s,t]} 
  b^\Lambda_n(i,k)^2\Bigr]
\nonumber \\&\le& 
\E\Bigl[\sum_{k\in\bar\Lambda}\sum_{n:T_n(k)\in [\bar s,t]} 
  b^{\bar\Lambda}_n(i,k)^2\Bigr]
\;=\; \E\,[\eta^{\bar\Lambda}_{[\bar s,t]}(i)]^2
\end{eqnarray}
where the inequality follows from the fact that the probabilities absorbed at
$\Lambda$ are dominated by the ones absorbed at $\bar\Lambda$: if
$\Lambda\subset\bar\Lambda$, then $b^\Lambda_n(i,k)\le b^{\bar\Lambda}_n(i,k)$.
This shows monotonicity in $\Lambda$ for \reff{F.19b}. Variances of martingales
are non decreasing in time, showing time monotonicity in \reff{F.19b} and
\reff{x20}. $\sqr$

\begin{corol}
  \label{xx11} 
There exist constants $C(i)$ such that for all $\Lambda$ and $s\le t$\\
(a) For $d\geq 3$,
  $\E\,[\eta^{\Lambda}_{[s,t]}(i)]^2 <C(i)$.\\
(b) For $d\ge 1$, $\E\,[\eta^{\Z^d}_{[s,t]}(i)
  - \eta^{\Z^d}_{[s,t]}(0)]^2<C(i)$.\\
(c) Assuming $p(0,0)=0$, for $d\ge 1$,
  $\E\,[\eta^{\Lambda\setminus\{0\}}_{[s,t]}(i)]^2<C(i)$.
\end{corol}

\proof (a) follows from \reff{F.19b}, \reff{F.19a} and \reff{bb2}. Obtain (b)
from \reff{x20}, \reff{F.19} and \reff{bb3} and (c) from \reff{F.19b} and \reff{x21x}.
$\sqr$

\section{Time and space convergence} \label{s4} The process
$(\eta^\Lambda_{[s,t]}\,:\, t\ge s)$ has ``flat boundary conditions'' outside
$\Lambda$ and ``flat initial condition'' at time $s$. We state the results for
this case and later comment about general boundary and initial conditions.  We
first show that under suitable conditions the process \reff{p20} is well defined
when $s=-\infty$ and it is in fact a stationary version of the harness process.
In particular, when the noise is Gaussian, the marginal law of this process at
any time $t$ has a Gibbs distribution with specifications \reff{p71} which are
also reversible for the harness processes with Gaussian noise. In one and two
dimensions there is no Gibbs measure with specifications \reff{p71} (see
\cite{G}, Chapter~13).  The harness process should not converge to a probability
measure in $d=1,2$ (delocalization); see Toom \cite{T} for the discrete-time
version.  However both the harness process pinned at the origin and the process
``as seen from the height at the origin'' converge to the pinned Gibbs measure
$\mu^{\Z^d\setminus\{0\}}$. The $L^2$ time convergence in $d\ge
3$ was proven by Hsiao \cite{Hs1}; we obtain the convergence bounds \reff{x3b}.

\begin{thm}
  \label{x1}
The following hold

{\bf A.s.\/ time convergence} Assume either (a) $d\ge 3$ or (b)
$\Lambda\neq\Z^d$ and $p(0,0)=0$. For each $t\in\R$,
$i\in\R^d$, as $s\to\infty$, $\eta^{\Lambda}_{[t-s,t]}(i)$ converges almost
surely to a random variable $\eta^{\Lambda}_t(i)$:
\begin{equation}
  \label{x2}
  \lim_{s\to\infty} \eta^{\Lambda}_{[t-s,t]}(i) \,=\,
  \eta^{\Lambda}_t(i) \quad\hbox{a.s.}
\end{equation}
For $d\ge 1$,
$\eta^{\Z^d}_{[t-s,t]}(i)-\eta^{\Z^d}_{[t-s,t]}(0)$ converges almost surely to a
random variable $\Delta^{\Z^d}_t(i)$:
\begin{equation}
  \label{x2a}
  \lim_{s\to\infty} \big[\eta^{\Z^d}_{[t-s,t]}(i)-\eta^{\Z^d}_{[t-s,t]}(0)\big] \,=\,
  \Delta^{\Z^d}_t(i) \quad\hbox{a.s.}
\end{equation}

{\bf $L_2$ time convergence} There exist positive constants $C, C(i)<\infty$
such that for $d\ge 1$, $\Lambda\subset\Z^d$ and $s\ge 0$,
\begin{eqnarray}
  \label{x3}
  \E(\eta^{\Lambda}_{[t-s,t]}(i) - \eta^{\Lambda}_{t}(i))^2\;\le\;
C\, s^{1-d/2}
\label{x3b}
\end{eqnarray}
(These bounds are relevant only for $d\ge 3$.)
\begin{equation}
  \label{x3c}
  \lim_{s\to\infty}\E(\eta^{\Lambda\setminus\{0\}}_{[t-s,t]}(i)-
  \eta^{\Lambda\setminus\{0\}}_{t}(i))^2 \;=\;0
\end{equation}
\begin{equation}
  \label{x3a}
  \E(\eta^{\Z^d}_{[t-s,t]}(i)-\eta^{\Z^d}_{[t-s,t]}(0) 
  - \Delta^{\Z^d}_{t}(i))^2
  \; \le \;
  C(i)\,s^{-d/2}
\end{equation}

{\bf Stationarity}
The processes $(\eta^\Lambda_t,\,t\in\R)$ and $(\Delta^{\Z^d}_t,\,t\in\R)$ are
stationary Markov with generators $L^\Lambda$ and $\widetilde L$ respectively,
where $\widetilde L$ is given later in \reff{q100}.

{\bf Uniqueness for finite $\Lambda$} If $\Lambda$ has a finite
number of points, then the law of $\eta^\Lambda_t$ is the unique invariant
measure for the process with generator $L^\Lambda$.

{\bf  $L_2$  space convergence} 
For either $d\ge 3$ or $\Lambda\neq\Z^d$,
\begin{equation}
  \label{x5}
  \lim_{\Lambda'\nearrow\Lambda} 
  \E(\eta^{\Lambda'}_t(i) - \eta^{\Lambda}_t(i))^2\;=\;0
\end{equation}
\end{thm}

\proof

{\bf A.s. time convergence} Fix $t\in\R$.  By Proposition~\ref{p777},
the process $(\eta^{\Lambda}_{[t-s,t]}(i),$ $s\ge t)$ is a martingale.  By
Corollary \ref{xx11} its variances are uniformly bounded under the given
conditions ---since the origin plays no special role, it is not a loss of
generality to assume that $0\notin\Lambda$). Analogously, the process as seen
from the height at the origin
$(\eta^{\Lambda}_{[t-s,t]}(i)-\eta^{\Lambda}_{[t-s,t]}(0),\,s\ge 0)$ is a
martingale with uniformly bounded variances under the given conditions.
Martingales with uniformly bounded variances converge almost surely \cite{HH}.

{\bf $L_2$ time convergence}
\begin{eqnarray}
  \lefteqn{\E(\eta^{\Lambda}_{[t-s,t]}(i) - \eta^{\Lambda}_{t}(i))^2}
  \\ &=& 
   \E(\eta^{\Lambda}_{[t-s,t]}(i))^2\, +\, \E(\eta^{\Lambda}_{t}(i))^2\,
  - \,2 
  \E(\eta^{\Lambda}_{[t-s,t]}(i)\eta^{\Lambda}_{t}(i)) \label{x28}\\ 
  &=& \E\,\sum_{k \in \Z^d}  
  \Bigl(\sum_{n: T_n(k) \in [t-s,t]}\,+\, \sum_{n: T_n(k) 
    \in (-\infty,t]}\,-\,\, 2 \sum_{n: T_n(k) \in [t-s,t]}\Bigr)\,
  b^\Lambda_{[t,T_n]}(i,k)^2\nonumber\\ 
  &=& \E\, \sum_{k \in \Z^d} \sum_{n: T_n(k)\in 
    (-\infty,t-s)}b^\Lambda_{[t,T_n]}(i,k)^2 \label{x66}
\;\le\;\E\, \sum_{k \in \Z^d} \sum_{n: T_n(k)\in
    (-\infty,t-s)}b^{\Z^d}_{[t,T_n]}(i,k)^2 \nonumber\\
  &=& \int_s^\infty  \,\P(D^0_u= 0)\,du \;<\; C s^{1-d/2}
\label{x27}  
\end{eqnarray}
where the second identity comes from \reff{F.1900}, the inequality from
\reff{F.19b} and the final identity can be shown as \reff{F.19a}. This shows
the inequality in \reff{x3b}.  

By the martingale property, 
\[
\E(\eta^{\Lambda\setminus\{0\}}_{[t-s,t]}(i) - 
\eta^{\Lambda\setminus\{0\}}_{t}(i))^2 
\;=\; \E(\eta^{\Lambda\setminus\{0\}}_{[t-s,t]}(i))^2 - 
\E(\eta^{\Lambda\setminus\{0\}}_{t}(i))^2
\]
which converges to 0 as $s\to\infty$ because it is an increasing bounded
sequence by Lemmas \ref{x70} and \ref{x11}. This shows
\reff{x3c}.
Analogously, using \reff{F.19} and \reff{bb3},
\begin{eqnarray}
  \label{x68}
\lefteqn{
\E(\eta^{\Z^d}_{[t-s,t]}(i)-\eta^{\Z^d}_{[t-s,t]}(0) 
  - \Delta^{\Z^d}_{t}(i))^2
}\nonumber\\
&=&2\,\int_{s}^\infty  \,(\P(D^0_u= 0)-\P(D^{i}_u= 0))\,du \;<\; C(i)
\,s^{-d/2} \label{x701} 
\end{eqnarray}

{\bf Stationarity} The construction of $\eta^\Lambda_t$ commutes with the
time-translation operator: $\eta^\Lambda_t(\omega+u) =
\eta^\Lambda_{t+u}(\omega) $, where $\omega=( (T_n(i), \eps_n(i), U_n(i))$:
$i\in \Z^d,\,n\in \Z)$ and $\omega+u :=( (T_n(i)+u, \eps_n(i), U_n(i))\,:\,
i\in \Z^d,\,n\in \Z)$ are identically distributed. The Markov property follows
as in \reff{I.13}.
\paragraph*{\bf Uniqueness} Let $\xi$ be a random configuration in $\R^{\L}$.
with invariant distribution for the process, then $\xi$ has the same law as
the random configuration
\begin{equation}
    \label{p48}
     \xi^\Lambda_{[s,t]} (i) := \sum_{j\in\Lambda}\;
  \sum_{n:T_n(j)\in [s,t)}
  \eps_n(j)\,b^\Lambda_n(i,j)
  +\sum_{j\in\Lambda}\;b_{[t,s]}^\Lambda(i,j)\xi(j)    
\end{equation}
(recall $b_{[t,s]}^\Lambda(i,j)= \P(B^{\Lambda,i}_{[t,s]} = j\vert\TT)$). Since
$\Lambda$ is finite and the walk $B^{\Lambda,i}_{[t,s]}$ is absorbed at
$\Lambda^c$, $b_{[t,s]}^\Lambda(i,j)$ goes to zero a.s.\/ as $s\to-\infty$ and so
does the second sum in \reff{p48}. This implies that $\xi$ and
$\eta^\Lambda_t$ (which is the limit of the first sum) have the same law.

{\bf $L_2$ space convergence}
Fix $i\in\Lambda$. Using \reff{F.1900} we get for $\Lambda\supset\Lambda'\ni i$,
\begin{eqnarray}
  \label{p62}
 \E\, (\eta^{\Lambda'}_t (i) - \eta^{\Lambda}_t(i))^2 
 &=& \E\,\sum_{j\in\Lambda'}\;
  \sum_{n:T_n(j)\le t}\,[b^\Lambda_n(i,j)\,-\, b^{\Lambda'}_n(i,j)]^2 
\end{eqnarray}
The summand in \reff{p62} is bounded by
$(b^{\Lambda}_n(i,j))^2+(b^{\Lambda'}_n(i,j))^2\le 2(b^{\Z^d}_n(i,j))^2$ which
is integrable in $d\ge 3$ by \reff{bb2} or if $\Lambda\neq \Z^d$ in $d=1,2$ by
\reff{x21x}. Then, since $ \lim_{\Lambda'\nearrow\Z^d}b^{\Lambda'}_n(i,j) =
b^{\Lambda}_n(i,j)$ a.s., \reff{p62} goes to zero as $\Lambda'\nearrow\Lambda$.
$\sqr$

\paragraph*{\bf The pinned process and the processes as seen from the
height at the origin}
The height at the origin of the process $\eta^{\Lambda\setminus\{0\}}_t(0)$
remains always equal to zero. For this reason, we call it the process
\emph{pinned at zero}. 

For fixed $s$, the process $(\eta^{\Z^d}_{[s,t]}- \eta^{\Z^d}_{[s,t]}(0),\,t\ge
s)$ is called the \emph{process as seen from the height at the origin}. Its
generator is
\begin{eqnarray}
  \label{q100}
  \widetilde L f(\eta) &=& \sum_{i\neq 0} \int G(d\eps)
  [f(P_i(\eta)+\sigma\eps e_i) - f(\eta)] \nonumber\\
&&\qquad+\;
  \int G(d\eps)
  \Bigl[f\Bigl(\eta-\Bigl(\sum_{\ell\neq 0} p(0,\ell) \eta(\ell)+\sigma\eps
  \Bigr) \sum_{j\neq
  0}e_j\Bigr) - f(\eta)\Bigr] 
\end{eqnarray}
The first term corresponds to updatings of sites other than the origin while the
second one corresponds to the shift all sites suffer when the origin is updated.

\paragraph*{\bf Convergence to the invariant measure} 
Due to the time stationarity of the marked Poisson processes, the law of
$\eta^{\Lambda}_{[s,t]}$ depends only on $t-s$, and in particular for
\emph{each} $t\ge 0$, $\eta^{\Lambda}_{[-t,0]}$ has the same law as
$\eta^{\Lambda}_{[0,t]}$. Hence, for cylinder Lipschitz functions $f$ for
which there exists a finite positive $\alpha$ satisfying
$|f(\eta)-f(\eta')|\le \alpha \big(\sum_k(\eta(k)-\eta'(k))^2\big)^{1/2}$
depending on the coordinates in the finite set Supp$(f)\subset \Z^d$,
\begin{eqnarray}
   \label{x7}
|\E f(\eta^{\Lambda}_{[0,t]}) - \mu^\Lambda f|
&=&  \big|\E \big(f(\eta^{\Lambda}_{[-t,0]}) -
f(\eta^{\Lambda}_{[-\infty,0]})\big)\big|\nonumber\\
&\le&  \E  \Bigl(\alpha\sum_{i\in\rm{Supp}(f)} 
(\eta^{\Lambda}_{[-t,0]}(i) -
\eta^{\Lambda}_{[-\infty,0]}(i))^2\Bigr)^{1/2}\nonumber\\
&\le& \Bigl(\alpha \sum_{i\in\rm{Supp}(f)} 
\E (\eta^{\Lambda}_{[-t,0]}(i) -
\eta^{\Lambda}_{[-\infty,0]}(i))^2\Bigr)^{1/2}\nonumber\\
&\le& \big(|\hbox{Supp}(f)|\,\alpha    
\,C\,t^{-1+d/2}\big)^{1/2} \nonumber
\end{eqnarray}
by \reff{x3b}. The last bound is relevant only in $d\ge 3$. Analogously, using
\reff{x3a},
\begin{eqnarray}
   \label{x7a}
|\E f(\eta^{\Z^d}_{[0,t]}-\eta^{\Z^d}_{[0,t]}(0)) - \mu^{\Z^d\setminus\{0\}} f|
&\le& \big(|\hbox{Supp}(f)|\,\alpha 
\,C\,t^{-d/2} \big)^{1/2}\nonumber
\end{eqnarray}

\paragraph*{\bf Other initial and boundary conditions}
Let $\Lambda\subset\Z^d$ and
\begin{equation}
  \label{o1}
  \Gamma_\Lambda := \{\gamma\,:\,\sum \bar b_\Lambda(i,j) \gamma(j)
  <\infty,\, \hbox{ for all }i\in\Lambda\} 
\end{equation}
where $b_\Lambda(i,j)$ is the probability that a continuous time
random walk, with rates $p$, absorbed at the sites of $\Lambda^c$,
starting at $i\in\Lambda$ is absorbed at site $j\in\Lambda^c$.

Let $\gamma\in\Gamma_\Lambda$ and $\zeta\in\Xi_\Lambda$ given in
\reff{o3}. Due to the linear property of the dynamics, the process
$\eta^{\Lambda,\gamma,\zeta}_{[s,t]}$ with initial configuration
$\eta^{\Lambda,\gamma,\zeta}_{[s,s]}=\zeta$ at time $s$ and boundary
conditions $\gamma$ can be seen as the sum of a process with flat
boundary and initial conditions plus a ``no noise'' harness process.

The process
$(\overline\eta^{\Lambda,\gamma,\zeta}_{[s,t]}\,:\,t\ge s)$ with
initial configuration $\overline\eta^{\Lambda,\gamma,\zeta}_{[s,s]}=\zeta$ at
time $s$ and generator
\begin{equation}
  \label{x35}
  \overline L^{\Lambda,\gamma} f(\eta) = \sum_{i\in\Lambda} 
  [f(P_i(\eta_\Lambda\gamma_{\Lambda^c})) - f(\eta)] 
\end{equation}
is called the \emph{no noise harness process}; it has $\gamma$
boundary conditions outside $\Lambda$. This is just a harness process
with noise distribution concentrating mass on the point 0 so that the
updating of site $i$ is done using only the $P_i$ average of the other
heights.  It is still a stochastic process because the updating times
are governed by the Poisson processes $\TT$.  Let
$\HH^{\Lambda,\gamma}$ be the set of harmonic functions for $p$ on
$\Lambda$ with $\gamma$ boundary conditions:
\[
  \HH^{\Lambda,\gamma}\,:=\, \Bigl\{h\in\R^{\Z^d}\,:\, \sum_j p(i,j)h(j) =
  h(i),\,  
  i\in\Lambda;\; h(j)=\gamma(j),\,i\in\Lambda^c 
  \Bigr\} \nonumber
\]
Measures concentrating mass on $\HH^{\Lambda,\gamma}$ are invariant for the
no-noise process $\overline\eta^{\Lambda,\gamma,\cdot}_{[s,t]}$. Some
questions naturally arise here: Do the invariant measures for the no-noise
process concentrate mass on $\HH^{\Lambda,\gamma}$? Does this process converge
to one of the invariant measures? If yes, what is the speed of convergence?

We have the following decomposition
\begin{equation}
  \label{x36}
  \eta^{\Lambda,\gamma,\zeta}_{[s,t]} 
\;=\; \eta^{\Lambda}_{[s,t]}+ \overline\eta^{\Lambda,\gamma,\zeta}_{[s,t]}
\end{equation}
Notice however that both processes use the same Poisson epochs. 
Measures in the set
\begin{equation}
  \label{x34}
  \II^{\Lambda,\gamma} = \{\hbox {law of }\eta^{\Lambda}_t + h \,:\,
  h\in\HH^{\Lambda,\gamma}\} 
\end{equation}
are invariant for the process $\eta^{\Lambda,\gamma}_{[0,t]}$. Are all
invariant measures convex combinations of the measures in
$\II^{\Lambda,\gamma}$? What are the domain of attraction of the measures in
$\II^{\Lambda,\gamma}$? 

\paragraph*{\bf Uniqueness} 
For $d\ge 3$, we conjecture that the law of $\eta^{\Z^d}_t$ is the unique
ergodic invariant measure with zero mean (that is, such that $\E\eta_t(i)=0$
for all $i$) for the process with generator~$L^{\Z^d}$. Hsiao \cite{Hs1} has
proven that the law of $\eta^{\Z^d}_t$ is the unique invariant measure with
zero mean and uniformly bounded second moment. For $d=1,2$, we conjecture that
the law of $\eta^{\Z^d\setminus\{0\}}_t$ is the unique ergodic (here we mean
for the height differences) invariant
measure with zero mean for the process with generator~$L^{\Z^d\setminus\{0\}}$
and the unique ergodic measure with
mean zero invariant for the pinned process $\eta^{\Z^d}_t-\eta^{\Z^d}_t(0)$.

\paragraph*{\bf A.s.\/{} space convergence}
Let $(\Lambda_m\,:\,m\ge 0)$ be an increasing family of sets such that
$\Lambda_m\nearrow\Lambda$. Assuming as extra condition that $G$ is Gaussian,
we exhibit a family of random configurations $(\xi^{\Lambda_m}_t\,:\,m\ge 0)$
with marginal laws $\mu^{\Lambda_m}$ converging almost surely as $\Lambda_m$
increases to $\Lambda$. As noted by the referee, the \emph{existence} of such
a sequence is guaranteed by the Skorohod representation theorem; our aim here
is to explicitely construct it.

Fix $\Lambda_m$ and the Poisson configuration $\TT$ and call
$b^m_n(i,j):=b^{\Lambda_m}_n(i,j)$ (this is a function of $\TT$). By \reff{p20},
$\eta^m_{[s,t]}(i):=\eta^{\Lambda_m}_{[s,t]}(i)$ is a sum of the independent
Gaussian random variables $\eps_n(j)\,b^m_n(i,j)$, for $n$ such that $T_n(j)\le
t$ and $j\in\Z^d$. 

Since $b^m_n(i,j)$ is non decreasing in $m$ we can define $a^0_n(i,j)=0$ and
for $m\ge 1$,
\begin{equation}
  \label{x80}
  a^m_n(i,j) 
  := \Bigl(b^{m}_n(i,j)^2-b^{m-1}_n(i,j)^2\Bigr)^{1/2}
\end{equation}
(so that $\sum_{\ell=1}^m (a^\ell)^2 = (b^m)^2$).  Let $Z_{n}^\ell(j)$ be a
sequence of independent and identically distributed centered Gaussian random
variables of variance 1 and let 
\begin{equation}
  \label{x81}
  W^m_n(i,j):= \sum_{\ell=1}^m a^\ell_n(i,j) Z_{n}^\ell(j)\,.
\end{equation}
Hence $W^m_n(i,j)$ are independent Gaussian random variables,
\begin{equation}
  \label{x82}
   W^m_n(i,j) \stackrel d= \eps_n(j) b^{m}_n(i,j)
\end{equation}
and the random configuration $\xi^m_t$ defined by
\begin{equation}
  \label{x83}
\xi^m_t(i):= \sum_j\sum_n  W^m_n(i,j) 
\end{equation}
has the same law as $\eta^m_t$. 

\begin{propos} 
\label{xx4} 
Assume $G(dx) = (2\pi)^{-1/2}e^{-x^2/2}dx$ (Gaussian noise). Then
 for either $d\ge 3$ or
$\Lambda\neq\Z^d$,
\begin{equation}
  \label{x4}
  \lim_{m\to\infty} \xi^m(i) \,=\,
  \xi^{\Lambda}_t(i) \quad\hbox{a.s.}
\end{equation}
and in $d\ge 1$, for any $\Lambda$,
\begin{equation}
  \label{x4a}
  \lim_{m\to\infty} (\xi^m_t(i)-\xi^m_t(0)) \,=\,
  \xi^{\Lambda}_t(i)-\xi^{\Lambda}_t(0) \quad\hbox{a.s.}
\end{equation}
\end{propos}

\proof 
By Lemma \ref{x84} below,
$(\xi^m_t(i)\,,\,m\ge 1)$ is a martingale. Since it has uniformly bounded
second moments by Corollary \ref{xx11}, it converges almost surely. $\sqr$

\begin{lemma}
  \label{x84}
  For each $i\in\Z^d$, the family $(\xi^m_t(i)\,,\,m\ge 1)$ is a martingale
  for the filtration $\FF_m$ generated by the family of variables $\{T_n(j),\,
  \, (Z^\ell_n(j)\,,\ell\le m)\,;\,j\in\Lambda_m,\,T_n(j)\le t \}$.
\end{lemma}

\proof Take $m'\ge m$. Then
\begin{equation}
  \label{x85}
  \xi^{m'}_t(i)-\xi^{m}_t(i)
  :=
  \sum_{j\in\Lambda_{m'}}\sum_{\ell=m+1}^{m'} \sum_n
 a^\ell_n(i,j) Z_{n}^\ell(j)
\end{equation}
which conditioned to $\FF_m$ has mean zero because it is a weighted sum of
$Z_{n}^\ell(j)$'s that are independent of the weights and of those
$Z_{n}^\ell(j)$'s generating $\FF_m$. $\sqr$

\section{Reversibility and Gibbs measures}\label{s5}
Most results of the previous sections hold for any variance-1 noise and for
any finite range matrix $p$. With this generality the properties of the law of
$\eta^\Lambda_t$ (which is an invariant measure for the process) are not well
understood besides the knowledge of the covariances.  However, if we assume
\begin{equation}
  \label{z1}
  G(dx) = (2\pi)^{-1/2}e^{-x^2/2}dx \hbox{ (Gaussian noise) and }p(0,0) = 0,  
\end{equation}
then for finite $\Lambda$ the law of $\eta^\Lambda_t$ is the finite volume
Gibbs measure $\mu^\Lambda$ given by \reff{p71} and it is reversible for
$L^\Lambda$. These properties extend to infinite $\Lambda$ as well. This is
the contents of our next result.

\begin{thm}
\label{xx3} 
Assume \reff{z1}.
Then,

1) For either $d\ge 3$ or $\Lambda\neq\Z^d$, the distribution of
$\eta^\Lambda_t$ is the Gibbs measure $\mu^\Lambda$ with
specifications \reff{p71} and boundary conditions $\gamma\equiv 0$ and
the process $(\eta^\Lambda_t,\,t\in\R)$ is reversible.

2) In $d\ge 1$, the marginal (invariant) distribution of $\eta^{\Z^d}_t-
\eta^{\Z^d}_t(0)$ is the Gibbs measure $\mu^{{\Z^d}\setminus\{0\}}$ with
specifications \reff{p71} and $\gamma\equiv 0$ and the process
$(\eta^{\Z^d}_t- \eta^{\Z^d}_t(0),\,t\in\R)$ is reversible.
\end{thm}

The case $d\ge 3$ and $\Lambda=\Z^d$ is already contained in Hsiao \cite{Hs1}.

\proof 1) For finite $\Lambda$ the statements are proven in Lemma \ref{p102}
below. For infinite $\Lambda$ the existence of the infinite volume measure
$\mu^\Lambda$ with specifications \reff{p71} is proven by Spitzer \cite{Sp2};
alternatively it follows either from the $L^2$ space convergence \reff{x5} in
Theorem \ref{x1} or the \emph{a.s.} space convergence of Theorem~\ref{xx4}.
The reversibility of the limiting measure $\mu^\Lambda$ follows then as in
Lemma \ref{p102}.

2) The existence of the infinite volume Gibbs measure
$\mu^{{\Z^d}\setminus\{0\}}$ is proven by Spitzer \cite{Sp2}, see also Caputo
\cite{C}.  We do not have an alternative proof in this case. The reversibility
follows as in Lemmas \ref{p102} and \ref{x60} later.  $\sqr$

Spitzer \cite{Sp2} (see Caputo \cite{C} for the non nearest neighbor case) proved
that the covariances of $\mu^{\Lambda\setminus\{0\}}$ are given by
\begin{equation}
  \label{z2}
  \int \mu^{{\Lambda}\setminus\{0\}}(d\xi)\, \xi(i)\xi(j) 
  = \sum_{n\ge 0} \P(X^i_n=j,\tau^i>n)
\end{equation}
where $X^i_n$ is a random walk with probability transition matrix $p$ and
$\tau^i$ is the first time the walk hits the origin or $\Lambda^c$. This is
the expected number of visits to $j$ for the walk $X_n$ starting at $i$ before
being absorbed at $0$ or $\Lambda^c$. These covariances are finite in any
dimension: the number of visits to $j$ of the walk starting at $j$ is a
geometric random variable because after each visit the walk can be absorbed at
0 or (in dimensions $d\ge 3$) never visit $j$ again. In particular there
exist constants $C(i)$
\begin{equation}
  \label{l12}
  V^{\Lambda\setminus\{0\}}(i)=:
  \int\mu^{\Lambda\setminus\{0\}}(d\xi)\,\xi(i)^2<C(i)<\infty  
  \qquad \hbox{ for all }\Lambda\subset \Z^d 
\end{equation}

The next lemma is essentially contained in Theorem 3.3 of Hsiao
\cite{Hs1}. 
\begin{lemma}
\label{p102}
Assume \reff{z1} and $\Lambda$ finite. Then the Gibbs measure
$\mu^{\Lambda,\gamma}$ with Hamiltonian $H^{\Lambda}(\eta) =
\frac12\sum_{i,j} p(i,j) (\eta(i)-\eta(j))^2$ is reversible for each of the
generators
\begin{equation}
  \label{p100a}
  L^{\Lambda,\gamma}_k f(\eta) =  \int G(d\eps)
  [f(P_k(\eta_\Lambda\gamma_{\Lambda^c})+\eps e_k) - f(\eta)]\,, \qquad
  k\in\Lambda 
\end{equation}
\end{lemma}
\noindent
(For definitions of $P_k$ and $\eta_\Lambda\gamma_{\Lambda^c}$ see
\reff{p100}).

\proof Denote $\mu=\mu^{\Lambda,\gamma}$, $L_k=L^{\Lambda,\gamma}_k$ and
$\eta= \eta_\Lambda\gamma_{\Lambda^c}$. We need to show that $\mu(gL_kf) =
\mu(fL_kg)$ for any continuous bounded functions $f$ and $g$. By definition,
\begin{eqnarray}
  \label{p1i}
\int \mu(d\eta) g(\eta) L_kf(\eta) &=&\int \mu(d\eta) g(\eta)  
\int \frac{e^{-x^2/2}}{\sqrt{2\pi}}dx\, [f(P_k\eta+e_k x)
\,-\, f(\eta)]\nonumber\\
&=& \int \mu(d\eta)  
\int \frac{e^{-x^2/2}}{\sqrt{2\pi}}dx \,g(\eta) \,f(P_k\eta+e_k x)\,
- \,\mu(gf)
\end{eqnarray}
Let $\bar\eta(k):= \sum_{i\neq k} p(k,i) \eta(i)$
(this does not depend on $\eta(k)$). Then, 
\begin{eqnarray}
  \label{p2i}
\sum_{i\neq k}
p(k,i)(\eta(k)-\eta(i))^2
&=&\sum_{i\neq k} p(k,i)(\eta(i)-\bar\eta(k))^2 +(\eta(k)-\bar\eta(k))^2\nonumber
\end{eqnarray}
Hence,
\begin{eqnarray}
  \label{p3i}
\lefteqn{ \int \mu(d\eta)\,g(\eta) \,\int G(dx)
\,f(P_k\eta+e_k x)}\label{p6i}\\ 
&=&\int \prod_{\ell\neq k} d\eta(\ell)\nonumber\\ 
&&\quad\times\exp\Bigl(- \frac12\sum_{i,j\neq k}p(i,j)(\eta(j)-\eta(i))^2\;
-\;\frac12\sum_{i\neq k} p(k,i)(\eta(i)-\bar\eta(k))^2\Bigr)\nonumber\\ 
&&\qquad\times 
\int e^{- (\eta(k)-\bar\eta(k))^2/2}\,d\eta(k)
 \int \frac{e^{-x^2/2}}{\sqrt{2\pi}}dx \,
g(\eta)f(P_k\eta+e_k x)\label{p7i}
\end{eqnarray}
Change variables: $\eta' = P_k\eta+e_k x$ and $z=\eta(k)-\bar\eta(k)$. Since
$\eta'(i)=\eta(i)$ for $i\neq k$, the second line in \reff{p2i} remains
unchanged when substituting $\eta$ for $\eta'$. Noticing that
$x=\eta'(k)-\overline{\eta'}(k)$ and $\eta=P_k\eta'-e_k z$, \reff{p7i} reads
\begin{eqnarray}
  \label{p4i}
 &&\int \prod_{\ell\neq k} d\eta'(\ell) \nonumber\\
&&\quad\times\exp\Bigl(-  \frac12\sum_{i,j\neq k}p(i,j)(\eta'(j)-\eta'(i))^2
\;-\; \frac12\sum_{i\neq k}
p(k,i)(\eta'(i)-\overline{\eta'}(k))^2\Bigr)\nonumber\\  
&&\qquad\qquad\times 
\int \frac{e^{-z^2/2}}{\sqrt{2\pi}}\,dz 
\int e^{-(\eta'(k)-\overline{\eta'}(k))^2/2}dx \,
g(P_k\eta'-e_k z)f(\eta')\nonumber\\
&&=\; \int \mu(d\eta) \,f(\eta)\,\int G(dz) 
\,g(P_k\eta+e_k z)\,.
\end{eqnarray}
Subtracting $\mu(fg)$ in \reff{p6i} and \reff{p4i} we obtain $\mu(gLf) =
\mu(fLg)$. $\sqr$

\paragraph*{\bf Free boundary conditions}
To find the infinite volume measure $\mu^{{\Z^d}\setminus\{0\}}$ we need to
introduce a family of processes and measures with free boundary conditions.
Let $\Lambda\subset \Z^d$ and define
\begin{equation}
  \label{x32}
  \widetilde p^{\Lambda}(i,j) := \frac
  {p(i,j)}{\sum_{k\in\Lambda}p(i,k)}\,,\qquad i,j\in\Lambda
\end{equation}
that is, a transition matrix for a walk that remains in $\Lambda$. Let
$\wL^{\Lambda}$, $\wH^{\Lambda}$, $\wmu^{\Lambda}$ be the generator,
Hamiltonian and Gibbs measure defined with $\wp^{\Lambda}$. In the dynamics
defined by $\wL^{\Lambda}$ the mean is taken only inside $\Lambda$ (no
boundary conditions matter). Let also 
\begin{equation}
  \label{p101a}
  \wL^{\Lambda}_k f(\eta) :=  \int G(d\eps)
  [f(\wP_k(\eta_\Lambda)+\sigma\eps e_k) - f(\eta_\Lambda)]\,, \qquad
  k\in\Lambda \setminus\{0\}
\end{equation}
the one-site generator of site $k\in\Lambda \setminus\{0\}$. We are interested
in two processes: the process with free boundary conditions pinned at zero and
the process with free boundary conditions as seen from the height at the
origin. The former one has generator $\sum_{k\in\Lambda \setminus\{0\}}
\wL^{\Lambda}_k$, while the second has generator $\sum_{k\in\Lambda
  \setminus\{0\}} \wL^{\Lambda}_k+ \wL^{\Lambda,0}f(\eta)$ where the shift
generator $\wL^{\Lambda,0}$ is defined by
\begin{equation}
  \label{x44}
  \wL^{\Lambda,0}_0f(\eta):= \int G(d\eps)
  [f(\eta_\Lambda - (\wP_0(\eta_\Lambda)(0)+\eps) \one) - f(\eta_\Lambda)] 
\end{equation}
for $f$ not depending on $\eta(0)$, where $\one$ is the configuration
$\one(i)\equiv 1$.

\begin{lemma}
  \label{x60}
  Assume \reff{z1} and $\Lambda$ finite. Then the Gibbs measure
  $\wmu^{\Lambda,0}$ is reversible for each of the generators
  $\wL^{\Lambda}_k$ (and hence for the process pinned at zero with free
  boundary conditions) and for the shift generator $\wL^{\Lambda,0}$ (and
  hence for the free process as seen from the height at the origin).
\end{lemma}

\proof The proof that the measure $\wmu^{\Lambda,0}$ is reversible for
$\wL^{\Lambda}_k$ for $k\neq 0$ goes as the proof of Lemma~\ref{p102}.

To show that $\wmu^{\Lambda,0}$ is reversible for $\wL^{\Lambda,0}$ take $g$
and $f$ not depending on the height at the origin and compute
\begin{eqnarray}
  \label{x37}
\lefteqn{ \int \wmu^{\Lambda,0}(d\eta)\, g(\eta)\, \int G(dx)
  \,f(\eta - (\bar\eta(0)+x) \one)}\\
 &=&\int \prod_{\ell\neq 0} d\eta(\ell)\, \\
&&\quad\times\exp\Bigl(-\frac12\sum_{i,j\neq 0}
\wp(i,j)(\eta(j)-\eta(i))^2-\frac12\sum_{k\neq 0}
\wp(0,k)\eta(k)^2\Bigr) \nonumber\\ 
&&\qquad\qquad\times\int \frac{e^{-x^2/2}}{\sqrt{2\pi}}dx\, 
g(\eta)\,f(\eta-(\bar\eta(0)+x)\one)
\label{x37a}
\end{eqnarray}
Change variables: $z=\bar\eta(0)$ and $\eta' = \eta-(\bar\eta(0)+x)\one$. Then
$\overline{\eta'}(0) = -x$ and
\begin{eqnarray}
  \label{x40}
 \sum_{k\neq 0}\wp(0,k)\eta(k)^2 + x^2 = \sum_{k\neq 0}\wp(0,k)\eta'(k)^2 + z^2
\end{eqnarray}
So that \reff{x37} equals
\begin{eqnarray}
  \label{x38}
 &=&\int \prod_{\ell\neq 0} d\eta'(\ell) 
\,\exp\Bigl(-\frac12\sum_{i,j\neq 0}
\wp(i,j)(\eta'(j)-\eta'(i))^2-\frac12\sum_{k\neq 0}
\wp(0,k)\eta'(k)^2\Bigr) \nonumber\\ 
&&\qquad\qquad\times\int \frac{e^{-z^2/2}}{\sqrt{2\pi}}dz 
\,g(\eta'- (\overline{\eta'}(0)+z)\one)\,f(\eta')
\nonumber\\
&=&\int \wmu^{\Lambda,0}(d\eta)\, f(\eta)\, \int G(dz)
  \,g(\eta - (\bar\eta(0)+z) \one)\label{x41}
\end{eqnarray}
Subtracting $\wmu^{\Lambda,0}(fg)$ in \reff{x37} and \reff{x41} we obtain
$\wmu^{\Lambda,0} (g\wL^{\Lambda,0}
f) = \wmu^{\Lambda,0}(f\wL^{\Lambda,0} g) $. $\sqr$

\section*{Acknowledgements}
We thank Marina Vachkovskaia, Luiz Renato Fontes, Servet Mart\'{\i}nez and Yvan
Velenik for fruitful discussions. We also thank a referee for useful comments.

This paper is partially supported by FAPESP, CNPq, PRONEX.
BN is supported by FAPESP through grant No. 00/05134--5.

\vskip 5mm
\parskip 0pt
\obeylines
\parindent0pt
Pablo A.~Ferrari, Beat M.~Niederhauser,
IME USP,
Caixa Postal 66281, 
05311-970 - S\~{a}o Paulo,
BRAZIL 

Phones: +55 11 3091 6119, +55 11 3091 6129, Fax: +55 11 38144 135

{\tt pablo@ime.usp.br},
 
{\tt http://www.ime.usp.br/~pablo}

\end{document}